\documentclass{article}
\usepackage{graphicx}%
\usepackage{multirow}%
\usepackage{amsmath,amssymb,amsfonts}%
\usepackage{amsthm}%
\usepackage[title]{appendix}%
\usepackage{xcolor}%
\usepackage{textcomp}%
\usepackage{manyfoot}%
\usepackage{booktabs}%
\usepackage{algorithm}%
\usepackage{algorithmicx}%
\usepackage{algpseudocode}%
\usepackage{listings}%
\usepackage{fix-cm}
\usepackage[scr=rsfso]{mathalfa}
\usepackage{hyperref}
\usepackage{geometry}
\newtheorem{theorem}{Theorem}
\newtheorem{lemma}{Lemma}
\newtheorem{example}{Example}%
\newtheorem{remark}{Remark}%

\newtheorem{definition}{Definition}%

\title{Convergence of Descent Optimization Algorithms under Polyak-\L ojasiewicz-Kurdyka Conditions}
\author{G. C. Bento\footnote{Federal University of Goiás, Goiânia, Goiás, BR (glaydston@ufg)}  \and B. S. Mordukhovich\footnote{Wayne State University, Detroit, MI 48202, USA (aa1086@wayne.edu)} \and T. S. Mota\footnote{Federal University of Goiás, Goiânia, Goiás, BR (tiago.mota@discente.ufg)} \and Yu. Nesterov\footnote{Corvinus Centre for Operations Research, Corvinus Institute for Advanced Studies, Corvinus University of Budapest, and School of Data Sciences (Chinese University of Hong Kong (Shenzhen); part time).  Professor emeritus at CORE/INMA, 
UClouvain, Belgium.\\}}
\date{\today}

\begin{document}

\maketitle
\vspace*{-0.3in}

\begin{abstract} {This paper develops a comprehensive  convergence analysis for generic classes of descent algorithms in nonsmooth and nonconvex optimization under several conditions of the Polyak-\L ojasiewicz-Kurdyka (PLK) type. Along other results, we prove the finite termination of generic algorithms under the PLK conditions with lower exponents. Specifications are given to establish new convergence rates for inexact reduced gradient methods and some versions of the boosted algorithm in DC programming. It is revealed, e.g., that the lower exponent PLK conditions for a broad class of difference programs are incompatible with the gradient Lipschitz continuity for the plus function around a local minimizer. On the other hand, we show that the above inconsistency observation may fail if the Lipschitz continuity is replaced by merely the gradient continuity.
 \\[1ex]
}\vspace*{-0.1in}
\end{abstract}
\section{Introduction, Basic Definitions, and Discussions}\vspace*{-0.05in}

The paper examines two generic classes of descent algorithms in (generally nonsmooth and nonconvex) constrained optimization. Our main goals are to establish global convergence and derive convergence rates for such algorithms under appropriate conditions. It has been well recognized that in optimization theory and many applications that conditions of the Polyak-\L ojasiewicz-Kurdyka type play in crucial role in the study of global convergence and convergence rates of various numerical algorithms. 

 For functions $\phi$ of  class $C^{1,1}$ (i.e., continuously differentiable functions with Lipschitzian gradients), Polyak introduced in 1962 (the English translation of his paper was published in \cite{polyak1963gradient}) the condition 
\begin{equation*}
\|\nabla\phi(x)\|\ge(1/2M)|\phi(x)-\phi(\Bar{x})|^{1/2},\quad c>0,
\end{equation*}
and used it to prove a linear convergence of the gradient descent method in Hilbert spaces (Theorem~4 in  \cite{polyak1963gradient}). Independently, \L ojasiewicz \cite{Lojasiewicz1963} introduced the inequality
\begin{equation}\label{pl}
\|\nabla\phi(x)\|\geq b\,|\phi(x)-\phi(\Bar{x})|^q,\quad b:=1/M(1-q),\quad q\in[0,1),
\end{equation}
for real analytic functions in the finite-dimensional framework of semialgebraic geometry with no applications to optimization. The gradient inequality \eqref{pl} is referred to (especially in the literature on machine learning and computer science) as the {\em Polyak-\L ojasiewicz} (PL) condition; see, e.g., \cite{karimi2016linear}. The subsequent algebraic-geometric extension of \eqref{pl} was developed by Kurdyka \cite{Kurdyka1998} for the class of definable differentiable functions on $\mathbb{R}^n$. A nonsmooth extension of the latter condition was first proposed by Bolt\'e et al.\ \cite{bolte2007lojasiewicz} in the form of Definition~\ref{def:kl}(ii) under the name of {\em Kurdyka-\L ojasiewicz} (KL) inequality without any reference to the pioneering work by Polyak. We suggest using the name of {\em Polyak-\L ojasiewicz-Kurdyka} (PLK) conditions for the properties of this type formulated below.  More details on the employed subgradient constructions and their calculi can be found in \cite{mordukhovich2006variational,mordukhovich2018variational,rw}. We also refer the reader to \cite{li-pong} for
important calculus rules for PLK exponents and their applications.\vspace*{0.05in}

\begin{definition}[\bf PLK conditions]\label{def:kl} {\rm Let $\phi:\mathbb{R}^n \to\overline{\mathbb R}:=\mathbb{R}\cup\{\infty\}$ be an extended-real-valued lower semicontinuous (l.s.c.) function with the domain  ${\rm dom}\,\phi:=\{x\in\mathbb{R}^n\;|\;\phi(x)<\infty\}$. We say that the function $\phi$ satisfies:

{\bf(i)} The {\em basic Polyak-\L ojasiewicz-Kurdyka} (PLK) {\em condition} at $\Bar{x}\in{\rm dom}\,\phi$ if there exist a number $\eta \in (0,\infty)$, a neighborhood $U$ of $\Bar{x}$, and a concave continuous function $\varphi: [0,\eta] \to[0,\infty)$, called the {\em desingularizing function}, such that
\begin{equation}\label{kl0}
\varphi(0)=0,\quad\varphi\in C^1(0,\eta),\quad\varphi'(s)>0\;\mbox{ for all }\; s\in(0,\eta),\;\mbox{ and }\;
\end{equation}
\begin{equation}\label{desigualdadeklND}
\varphi'\big(\phi(x)-\phi(\Bar{x})\big)\mbox{dist}\big(0,\partial\phi(x)\big)\ge 1\;\mbox{ for all }\;x\in U\cap[\phi (\Bar x) <\phi(x)<\phi(\Bar{x})+\eta]
\end{equation}
where $\partial\phi(\Bar x)$ stands for the {\em Mordukhovich/limiting subdifferential} of  $\phi$ at $\Bar{x}$ defined by
\begin{equation}\label{ls}
\partial\phi(\Bar{x}):=\Big\{v\in\mathbb{R}^n\;\Big|\;\exists x_k\stackrel{\phi}{\to}\Bar{x}, \,v_k\in\widehat\partial\phi(x_k),\quad v_k\to v\;\mbox{ as }\;k\to\infty\Big\},
\end{equation}
with $x_k\stackrel{\phi}{\to}\Bar{x}$ meaning that $x_k\to \Bar{x}$, $\phi(x_k)\to\phi(x)$, and with $v_k\in\widehat\partial\phi(x_k)$ meaning that 
$$\limsup_{u\to x_k}\frac{\phi(u)-\phi(x_k)-\langle v_k,u-x_k\rangle}{\|u-x_k\|}\ge 0,\;k\in\mathbb N:=\big\{1,2,\ldots\}.
$$

{\bf(ii)} The {\em symmetric PLK condition} at $\Bar{x}$ if $\phi$ is continuous around $\Bar{x}$ and $\partial\phi$ is replaced in \eqref{desigualdadeklND} by the {\em symmetric subdifferential} of $\phi$ at $\Bar{x}$ defined by
\begin{equation}\label{sym}
\partial^0\phi(\Bar{x}):=\partial\phi(\Bar{x})\cup\big(-\partial(-\phi)(\Bar{x})\big).
\end{equation}

{\bf(iii)} The {\em strong PLK condition} if $\phi$ is Lipschitz continuous around $\Bar{x}$ and $\partial$ is replaced in \eqref{desigualdadeklND} by the {\em Clarke/convexified subdifferential} of $\phi$ at $\Bar{x}$ given by
\begin{equation}\label{cl}
\overline{\partial}\phi(\Bar{x}):={\rm co}\,\partial\phi(\Bar{x}),
\end{equation}
where ``{\rm co}" stands for the convex hull of the set in question.

{\bf(iv)} The {\em exponent} versions of the {\em PLK conditions} in {\rm(i)--(iii)} if the desingularizing function in \eqref{kl0}  and \eqref{desigualdadeklND} is selected in the form $\varphi(t)=M t^{1-q}$, where $M$ is a positive constant, and where $q\in[0,1)$. We refer to the case where $q\in(0,1/2)$ as the {\em PLK conditions with lower exponents}.}
\end{definition}\vspace*{0.05in}

In \cite{attouch2013convergence}, Attouch, Bolt\'e and Svaiter studied minimization  of l.s.c.\ functions $\phi\colon\mathbb{R}^n\to\overline{\mathbb R}$ by using a {\em generic} class of {\em descent methods} satisfying the following properties:
\begin{itemize}
\item[\textbf{($\mathcal{H}1$)}]\textit{Sufficient decrease}: for each $k \in \mathbb{N}$, we have
\begin{equation}\label{desigualdade1kNDH}
\phi(x^{k+1}) + a \|x^{k+1} - x^k\|^2 \le\phi(x^k).
\end{equation}
\item[\textbf{($\mathcal{H}2$})] \textit{Relative error}: for each $k\in \mathbb{N}$, there exists
\begin{equation}\label{fran}
w^{k+1} \in \partial\phi(x^{k+1})\;\mbox{ with }\;\|w^{k+1}\| \leq b \|x^{k+1} - x^k\|,
\end{equation}
\end{itemize}
with $a,b>0$. It is shown in \cite{attouch2013convergence} that a great variety of important algorithms of optimization satisfy conditions $({\cal H}1)$ and $({\cal H}2)$. Employing the basic PLK property from Definition~\ref{def:kl}(i), a qualitative convergence analysis, with arriving at the {\em limiting/$($M$)$ordukhovich-stationary point} as $0\in\partial\phi(\Bar{x})$, is developed in \cite{attouch2013convergence} for the general class of descent algorithms satisfying $({\cal H}1)$ and $({\cal H}2)$ while without establishing any convergence rate.\footnote{In fact, an additional technical and nonrestrictive assumption, labeled as  $({\cal H}4)$ in what follows, was imposed in \cite{attouch2013convergence} and the subsequent publications discussed below.} A detailed  qualitative (with convergence rates) analysis of the generic class of descent algorithms under the exponent PLK conditions married to $({\cal H}1)$ and $({\cal H}2)$ was developed in \cite{frankel2015splitting}, where the 
the obtained results were similar to \cite{attouch2009convergence} established for the proximal algorithm; namely, the finite step convergence if $q = 0$, linear convergence if $q\in (0,1/2]$, and polynomial convergence if $q \in (1/2,1)$.\vspace*{0.03in}

A major motivation for our study came from the paper by Arag\'on-Artacho and Vuong \cite{aragon2020boosted} who observed that the relative error condition $({\cal H}2)$ may fail for the  {\em Boosted Difference of Convex Algorithm} (BDCA) proposed in \cite{Aragon2018} for DC (difference of convex) programs. As shown in \cite{aragon2020boosted}, the convergence analysis of \cite{attouch2013convergence} can be given for BDCA {\em provided that} the relative error condition \eqref{fran} is replaced by 
\begin{equation}\label{fran1}
\mbox{there exists }\;w^{k}\in\overline{\partial}\phi(x^{k})\;\mbox{ with }\;\|w^{k}\|\leq b \|x^{k+1} - x^k\|,\quad k\in\mathbb N,
\end{equation}
expressed in terms of the convexified subdifferential \eqref{cl}. To furnish this, it is assumed in  \cite{aragon2020boosted} that the functions $g$ and $h$ in the DC decomposition $\phi=g-h$ are strongly convex with $g$ being of class ${\cal C}^{1,1}$, and that the basic PLK condition in \eqref{desigualdadeklND} is replaced by its strong version from Definition~\ref{def:kl}(iii). Under these requirements,  \cite{aragon2020boosted} establishes the convergence of BDCA iterates to a {\em $(C)$larke-stationary point} $0\in\overline{\partial}\phi(\Bar x)$.\vspace*{0.03in}

In the current paper, besides establishing new results on the convergence rate for generic algorithms satisfying conditions $({\cal H}1)$  and $({\cal H}2)$ under the PLK properties, we also investigate another generic class of algorithms where the relative error condition $({\cal H}2)$ is replaced by the following one:\\[1ex]
{\bf{($\mathcal{H}3$)}} {\em Modified relative error}: for each $k\in \mathbb{N}$, there exists
\begin{equation}\label{fran2}
w^{k}\in\partial\phi(x^{k})\;\mbox{ with }\;\|w^{k}\|\leq b \|x^{k+1} - x^k\|
\end{equation}
expressed in terms of the limiting subdifferential \eqref{ls}. The convergence analysis conducted below allows us to achieve, under the basic PLK property, the global convergence of any algorithm satisfying $({\cal H}1)$ and $({\cal H}3)$ to an $M$-stationary point $0\in\partial\phi(\Bar{x})$. This answers in the affirmative the question posted in \cite[Remark~4.5]{aragon2020boosted}, even with the improvement \eqref{fran2} of the stronger condition in \eqref{fran1} for general l.s.c.\ functions. 

While the global convergence results under $({\cal H}3)$ are parallel to those developed under $({\cal H}2)$, this is not the case when it comes to convergence rate analysis. In what follows, we prove that the {\em lower exponent PLK} condition $(0<q<1/2)$ always yields the {\em finite termination} of any descent algorithm satisfying $({\cal H}1)$ and $({\cal H}3)$. On the other hand, we present an example showing that the finite termination {\em fails} under PLK lower exponents for the proximal method, which surely satisfies $({\cal H}1)$ and $({\cal H}2)$.

Moreover, it comes as a surprise that in the case of difference programs to minimize $\phi=g-h$, where $h\colon\mathbb{R}^n\to\overline{\mathbb R}$ is convex while $g$ is of class ${\cal C}^{1,1}$ without any convexity assumption, there exists {\em no desingularizing function} of the type $\varphi(t)=t^{1-q}$ with $q\in(0,1/2)$ for which the PLK condition holds at {\em local minimizers} of $\phi$.\vspace*{0.05in}

The rest of the paper is organized as follows. Section~\ref{sec:conv} establishes {\em global convergence} of descent algorithms satisfying the conditions in $({\cal H}1)$ and $({\cal H}3)$, and in parallel to those in $({\cal H}1)$ and $({\cal H}2)$, to the corresponding stationary points defined by the subdifferentials in the PLK properties from Definition~\ref{def:kl}(i--iii).

In Section~\ref{sec:rel-error}, we develop a novel {\em rate convergence analysis} of the generic algorithms satisfying  $({\cal H}1)$ and $({\cal H}2)$ under the exponent PLK conditions. A striking result of this section shows that the {\em lower exponent} PLK case exhibits {\em either the finite termination, or the superlinear} global convergence rate. It is a serious improvement of the linear convergence rate known in such settings. The failure of finite termination is demonstrated here by a simple example for the proximal algorithm.

Section~\ref{sec:rel-modif} continues our rate convergence analysis under the exponent PLK conditions, while now with replacing $({\cal H}2)$ but its {\em modified relative error} counterpart $({\cal H}3)$. The main novelty here is to show that any algorithm satisfying $({\cal H}1)$ and $({\cal H}3)$ {\em always  terminates in finite steps} under the lower exponent PLK conditions. Moreover, we present an additional assumption under which the latter conditions are {\em inconsistent} with the properties in  $({\cal H}1)$ and $({\cal H}3)$.

Section~\ref{sec:specif} contains some specifications of applications to particular algorithms that lie within the generic class of algorithms satisfying the conditions in $({\cal H}1)$ and $({\cal H}3)$. These specific algorithms concern nonsmooth versions of the inexact reduced gradient methods recently proposed in \cite{Mordukhovich2023inexact} as well as the aforementioned BDCA. It will be demonstrated that the obtained results for the generic descent methods lead us to significant improvements and extensions of the previously known ones.

Section~\ref{sec:KL} addresses the fulfillment of the PLK conditions with lower exponents $q\in(0,1/2)$ in the case of minimizing functions that admit the difference decomposition $\phi=g-h$. We show that the ${\cal C}^{1,1}$ structure of $g$ and the convexity of $h$ is {\em inconsistent} with the lower exponent PLK of $\phi$ at local minimizers. Section~\ref{sec:final} summarizes the major achievements of the paper and discusses some directions of the future research. 

The notation of this paper is standard in variational analysis and optimization. Recall that $B(x,r)$ stands for the open ball centered at $x$ with radius $r>0$.\vspace*{-0.1in}

\section{Global Convergence of Generic Algorithms}\label{sec:conv}\vspace*{-0.05in}

This section studies in parallel the two classes of descent algorithms satisfying the generic conditions  in {\bf(a)}: $({\cal H}1)$ and $({\cal H}2)$, and in {\bf(b)}: $({\cal H}1)$ and $({\cal H}3)$, respectively. Similarly to \cite{attouch2013convergence}, we add to our convergence analysis the following technical assumption:\\[1ex]
$({\cal H}4)$ {\em Continuity condition}:
There are a subsequence $\{x^{k_j}\}$ and $\Bar{x}\in{\rm dom}\,\phi$ such that
\begin{equation*}
x^{k_j} \longrightarrow\Bar{x}\;\text{ and }\;\phi(x^{k_j}) \longrightarrow \phi(\Bar{x})\;\mbox{ as }\;j\longrightarrow\infty.
\end{equation*}
From now on, unless otherwise stated, condition $(\mathcal{H}4)$ is  assumed to be satisfied.\vspace*{0.05in}

The major theorem below tells us that the basic PLK property from Definition~\ref{def:kl} ensures the global convergence of generic algorithms in both cases (a) and (b). In Remark~\ref{rem:strong}, we discuss the corresponding modifications of the convergence results under the symmetric and strong PLK conditions formulated in Definition~\ref{def:kl}(ii,iii).\vspace*{0.05in}

\begin{theorem}[\bf global convergence of generic algorithms under the basic PLK]\label{theo1}
Let $\phi :\mathbb{R}^n \longrightarrow\overline{\mathbb{R}}$ be a proper l.s.c.\ function bounded from below, and let the sequence $\{x^k\}$ be constructed by the generic algorithms satisfying either $(\mathcal{H}1)$ and $(\mathcal{H}2)$, or $(\mathcal{H}1)$ and $(\mathcal{H}3)$ properties. If the basic PLK condition holds at some accumulation point $\Bar{x}\in{\rm dom}\,\phi$ of $\{x^k\}$, then we have that
\begin{equation}\label{sum1}
\sum_{k=0}^{\infty}\left\|x^k-x^{k+1}\right\|<\infty,
\end{equation}
and that $\{x^k\}$ converges to $\Bar{x}$ as $k\to\infty$. Moreover, $\Bar{x}$ is an $M$-stationary point of $\phi$.
\end{theorem}\vspace*{-0.05in}
\begin{proof}
The result under $(\mathcal{H}1)$ and $(\mathcal{H}2)$ is known from \cite[Theorem~2.9]{attouch2013convergence}, and thus we proceed with the proof in the case where $(\mathcal{H}1)$ and $(\mathcal{H}3)$ hold. Since $\bar{x}$ is an accumulation point of $\{x^k\}$, there exists a subsequence $\{x^{k_j}\}$ of $\{x^k\}$ converging to $\bar x$ as $j\to\infty$. By the boundedness from below of $\phi$ and the decreasing of $\{\phi(x^k)\}$, the continuity condition $(\mathcal{H}4)$ implies that $\{\phi(x^k)\}$ converges to $\phi(\bar{x})$ as $k\to\infty$ with $\phi(\bar{x}) < \phi(x^k)$, $k\in \mathbb{N}$. In particular, we get
\begin{equation*}
\phi(\bar{x}) < \phi(x^k) < \phi(\Bar{x}) +\eta\;\mbox{ for large }\;k\in\mathbb N\;\mbox{ and small }\eta >0.
\end{equation*}
Taking the latter into account, define the sequence $\{b_k \}$ by
$$
b_k:= \left\| x^k - \bar{x} \right\| +\big(a^{-1} (\phi(x^k) - \phi(\Bar{x}))  \big)^{\frac{1}{2}} +b a^{-1} \varphi\big(\phi(x^k) - \phi(\Bar{x})\big)
$$
and deduce from the continuity of $\varphi$ that the origin $0\in\mathbb{R}^n$ is an accumulation point of $\big\{b_k\big\}$. Thus there exists $k_0 := k_{j_0}$ such that the inequalities
\begin{equation}\label{converg 4.2}
\left\| x^{k_0} - \bar{x} \right\| +a^{-1}\big(\phi(x^{k_0}) - \phi(\Bar{x})\big)^{\frac{1}{2}} +b a^{-1}\varphi\big(\phi(x^{k_0}) - \phi(\Bar{x})\big) < \epsilon,
\end{equation}
\begin{equation}\label{converg 4.3}
\phi(\Bar{x}) < \phi(x^{k_0}) < \phi(\Bar{x}) + \eta
\end{equation}
are satisfied and yield in turn the inclusions
$$
x^{k_0} \in B(\Bar{x},\epsilon) \cap \big[\phi(\Bar{x}) < \phi < \phi(\Bar{x})+\eta\big] \subset U \cap \big[\phi(\Bar{x}) < \phi < \phi(\Bar{x}) + \eta\big].
$$
Moreover, conditions \eqref{desigualdade1kNDH} and \eqref{converg 4.3} ensure that
\begin{equation}\label{converg 4.4}
0\leq \phi(x^{k_0 +1}) - \phi(\Bar{x}) \leq \phi(x^{k_0}) - \phi(\bar{x}) < \eta,
\end{equation}
which implies therefore the inequality
\begin{equation*}
\varphi\big(\phi(x^{k_0}) - \phi(\bar{x})\big) - \varphi\big( \phi(x^{k_0+1}) - \phi(\bar{x})\big)\geq \varphi'\big(\phi(x^{k_0}) - \phi(\bar{x})\big)\big( \phi(x^{k_0}) -\phi(x^{k_0+1})\big)
\end{equation*}
by the concavity of $\varphi$. Combining \eqref{desigualdade1kNDH} and \eqref{fran2} with $k:=k_0$, we have
$$
\phi(x^{k_0}) - \phi(x^{k_0+1}) \geq b a^{-1}\|x^{k_0 +1} - x^{k_0}\|\cdot\|w^{k_0}\|.
$$
The last two inequalities above lead us to
\begin{equation}\label{converg 4.6}
\begin{array}{ll}
&\varphi\big(\phi(x^{k_0}) - \phi(\bar{x})\big) - \varphi\big( \phi(x^{k_0+1}) - \phi(\bar{x})\big)\\
&\ge\displaystyle\frac{b}{a}\varphi'\big(\phi(x^{k_0}) - \phi(\bar{x})\big) \left\| w^{k_0} \right\|\left\| x^{k_0} - x^{k_0+1}\right\|
\end{array}
\end{equation}
with $w^{k_0}\in \partial f(x^{k_0})$, and so the basic PLK property of $f$ at $\Bar{x}$ yields $0 \notin \partial f(x^{k_0})$ while ensuring together with \eqref{desigualdadeklND} that
$$
\varphi'\big(\phi(x^{k_0}) - \phi(\bar{x})\big)\left\|w^{k_0} \right\|\ge\varphi' \big(\phi(x^{k_0}) - \phi(\bar{x})\big)\text{dist}\big(0,\partial \phi(x^{k_0})\big)\geq 1.
$$
Combining the latter with \eqref{converg 4.6} justifies the estimates
\begin{equation}\label{converg 4.7}
b a^{-1}\big(\varphi(\phi(x^{k}) - \phi(\bar{x})) - \varphi (\phi(x^{k+1}) -\phi(\bar{x}))\big) \geq \left\| x^{k} - x^{k+1}\right\|,\quad k=k_0.
\end{equation}

Our next step is to verify the inclusion
\begin{equation}\label{converg 4.8}
x^k \in B(\Bar{x}, \epsilon)\;\mbox{ for all }\;k\geq k_0.
\end{equation}
To proceed by induction, observe that we have already proved that \eqref{converg 4.8} holds for $k=k_0$ and now aim at showing that  $x^{k_0 +1} \in B(\Bar{x}, \epsilon)$. It follows from \eqref{desigualdade1kNDH} and \eqref{converg 4.4} that
\begin{equation}\label{converg 4.9}
\left\| x^{k_0} - x^{k_0 +1} \right\| \le\big(a^{-1}\big(\phi(x^{k_0}) - \phi(x^{k_0 +1})\big)\big)^{\frac{1}{2}} < \big(
 a^{-1}\big(\phi(x^{k_0}) - \phi(\Bar{x})\big)\big)^{\frac{1}{2}}.
\end{equation}
The last estimate along with the triangle inequality gives us
$$
\left\| \bar{x} - x^{k_0 +1} \right\| \leq \left\|
\Bar{x} - x^{k_0} \right\| + \left\| x^{k_0} - x^{k_0 +1} \right\| \leq \left\| \Bar{x} - x^{k_0} \right\| + \big(
a^{-1}\big(\phi(x^{k_0}) - \phi(\Bar{x}))\big)^{\frac{1}{2}}.
$$
Combining this with \eqref{converg 4.2} implies that $x^{k_0 +1}\in B(\Bar{x}, \epsilon)$. Take $j>1$ and assume that \eqref{converg 4.8} holds for all $k= k_0+1,\ldots,k_0 +j-1$. In this case, \eqref{converg 4.7} is satisfied for $k= k_0+1,\ldots, k_0 +j-1$. Therefore,
\begin{equation*}
\sum_{i=1}^{j-1} \left\|x^{i+k_0} - x^{i+k_0 +1} \right\| \le b a^{-1}\big( \varphi (\phi(x^{k_0 +1}) - \phi(\bar{x})) - \varphi\big(\phi(x^{k_0 +j}) - \phi(\bar{x}))\big).
\end{equation*}
Using the triangle inequality again tells us that
$$
\left\| \bar{x} - x^{k_0 +j} \right\| \leq \sum_{i=1}^{j-1} \left\|x^{i+k_0} - x^{i+k_0 +1} \right\| + \left\| x^{k_0} - x^{k_0 +1} \right\| + \left\| \bar{x} - x^{k_0} \right\|.
$$
Now we combine the last two inequalities with \eqref{converg 4.9} to get the estimates
$$
\left\| \bar{x} - x^{k_0 +j} \right\| \le b a^{-1} \varphi\big(\phi(x^{k_0 +1}) - \phi(\bar{x})\big) +\big(a^{-1} (\phi(x^{k_0 +1}) - \phi(\bar{x}))\big)^{\frac{1}{2}} +\left\| \bar{x} - x^{k_0}\right\|,
$$
$$\;\;\;\;\;\;\;\;\;\;\;\;\;\;\;\;\leq  b^{-1}a\varphi\big(\phi(x^{k_0}) - \phi(\bar{x})\big) + \big(a^{-1} (\phi(x^{k_0 +1}) - \phi(\bar{x}))\big)^{\frac{1}{2}} + \left\| \bar{x} - x^{k_0} \right\|,$$
where the second one follows from \eqref{converg 4.4} and the fact that $\varphi$ is increasing. The latter inequality along with \eqref{converg 4.2} ensures that $x^{k_0 +j}\in B(\bar{x},\epsilon)$, which verifies \eqref{converg 4.8} and shows that \eqref{converg 4.7} holds for all $k \geq k_0$. Thus we arrive at
\begin{equation*}
\sum_{k=k_0}^{\infty} \left\|x^{k} - x^{k+1} \right\| \leq b a^{-1} \varphi\big(\phi(x^{k_0}) -\phi(\Bar{x})\big)\;\mbox{ whenever }\;k>k_0,
\end{equation*}
which clearly yields \eqref{sum1}. This tells us that $\{x^k\}$ is a Cauchy sequence and hence converges to $\Bar{x}$ since it is an accumulation point of $\{x^k\}$. It follows from \eqref{fran2} and \eqref{sum1} that the sequence $\{w^k\}$ converges to zero. By definition \eqref{ls} of the limiting subdifferential $\partial\phi$ and the convergence $\phi(x^k)\to\phi(\Bar{x})$ as $k\to\infty$ shown above, we get that $0\in \partial \phi(\bar{x})$, i.e., $\bar{x}$ is an $M$-stationary point of $\phi$, and the proof is complete.
\end{proof}\vspace*{-0.05in}

The following remark reveals global convergence properties of generic descent methods under the {\em symmetric} and {\em strong} PLK conditions.\vspace*{0.03in}

\begin{remark}[\bf convergence under symmetric and strong PLK conditions]\label{rem:strong} $\,$ \vspace*{-0.07in}{\rm
\begin{itemize}
\item[\bf(i)] Theorem~\ref{theo1} resolves in the affirmative the question posed in \cite[Remark~4.5]{aragon2020boosted} about the possibility to provide a counterpart of the convergence analysis in \cite{attouch2013convergence} under the modified error bound condition \eqref{fran2} in (${\cal H}3)$.

\item[\bf(ii)] It follows from the proof of Theorem~\ref{theo1} that our analysis works with replacing \eqref{fran2} by \eqref{fran1} under the strong PLK condition from Definition~\ref{def:kl}(iii) and also by using the symmetric subdifferential \eqref{sym} in \eqref{fran2} under the symmetric PLK property from Definition~\ref{def:kl}(ii). However, such replacements lead us to  the less informative {\em C-stationary} $0\in\overline{\partial}\phi(\bar x)$ and {\em symmetric stationary} points $0\in\partial^0\phi(\bar x)$, respectively. Similar conclusions concerning the strong and symmetric PLK conditions can be made for generic algorithms in the case of $(\mathcal{H}2)$. Moreover, the same observation applies to the {\em convergence rate results} presented in Section~\ref{sec:rel-error} and Section~\ref{sec:rel-modif} below.
\end{itemize}}
\end{remark}\vspace*{-0.2in}


\section{Convergence Rates under Relative Error Condition}\label{sec:rel-error}\vspace*{-0.05in}


In this section, we provide explicit convergence rates in terms of the PLK exponents of the desingularizing function $\varphi(t) = t^{1-q}$, $q \in [0,1)$, for generic algorithms that satisfy the conditions in $({\cal H}1)$ and $({\cal H}2)$. While the convergence rate results obtained under the exponent PLK property for $q=0$ and $q\in[1/2,1)$ are generic extensions of those given for particular algorithms (see, e.g., \cite{Aragon2018,aragon2023coderivative,aragon2020boosted,attouch2009convergence,li,li-pong,Mordukhovich2023inexact,mordukhovich2024} and the references therein), it is worth highlighting that the main challenge lies in the case where $q\in(0,1/2)$ for which we provide significant convergence rate improvements.

To proceed, we  first present two technical lemmas that will be instrumental for our subsequent results. The idea behind the proof of the first lemma can be found within the proof of \cite[Theorem~3.2]{Absil2005}. As for the second lemma, to the best of our knowledge, it has not been previously documented in the literature and will play a crucial role in our generic convergence rate analysis.\vspace*{0.05in}

\begin{lemma}[\bf estimates for monotone sequences]\label{propositionAbsil2005}
Let $\{a_k\}\subset[0,\infty)$ be a monotonically decreasing sequence, and let $q\in (0,1)\cup(1,2)$. Then we have the estimate
\begin{equation*}
\sum_{j=k}^{k+l} \frac{a_j - a_{j+1}}{a_j^q} \leq \frac{1}{1-q}\left(a_k^{1-q} - a_{k+l+1}^{1-q}\right)\;\mbox{ whenever }\;k,l \in \mathbb{N}.
\end{equation*}
\end{lemma}\vspace*{-0.1in}
\begin{proof}
Since $\{a_k\}$ is monotonically decreasing, it follows that $a_j \geq t$ for all $t \in [a_{j+1}, a_j]$ and $j\in\mathbb{N}$. This leads us to
\begin{equation*}
\frac{1}{(a_j)^q} \leq \frac{1}{t^q}, \qquad t \in [a_{j+1}, a_j], \quad q \in (0,1),
\end{equation*}
which implies in turn that
\begin{equation*}
\frac{a_j - a_{j+1}}{a_j^q} =\int_{a_{j+1}}^{a_j} \frac{1}{a_j^q}dt\leq \int_{a_{j+1}}^{a_j} \frac{1}{t^q}dt.
\end{equation*}
The latter readily ensures that 
\begin{equation*}
\frac{a_j - a_{j+1}}{a_j^q}\leq\frac{1}{1-q } \left(a_j^{1- q} - a_{j+1}^{1-q}\right), \quad j\in\mathbb{N},
\end{equation*}
and thus we arrive at the claimed estimate.
\end{proof}

The next lemma is essential to establish the main results of this section on either finite or superlinear convergence of generic algorithms satisfying  the conditions in $({\cal H}1)$ and $({\cal H}2)$ under the PLK lower exponents.\vspace*{0.05in}

\begin{lemma}[\bf convergence rates for monotone sequences]\label{lemmatec}
Let $\sigma>0$ and $p \in (0,1)$, and let $\{a_k\}\subset [0,\infty)$ be a monotonically decreasing sequence. Assume that
\begin{equation}\label{comparacaok+1}
a_{k+1} \leq a_k - \sigma a_{k+1}^{p} \quad k \in \mathbb{N}.
\end{equation}
If $\{a_k\}$ converges to $\bar{a}$, then $\bar{a}=0$ and one of the following two conditions holds:\vspace*{-0.05in}
\begin{itemize}
\item[\bf(i)] The sequence $\{a_k\}$ converges superlinearly to zero.
\item[\bf(ii)] There exist $k_0 \in \mathbb{N}$ such that $a_k = 0$ for all $k\geq k_0$.
\end{itemize}
\end{lemma}\vspace*{-0.1in}
\begin{proof}
Since $a_k\geq 0$ for all $k \in \mathbb{N}$, we have $\bar{a}\ge 0$ be the limit of this sequence. On the other hand, the passage to the limit in \eqref{comparacaok+1} with taking into account that $\sigma>0$ yields
$$
\bar{a}\leq \bar{a} - \sigma \bar{a}^p\Longrightarrow\sigma \bar{a}^p\le 0,
$$
which implies that $\bar{a} = 0$. Fix $p\in (0,1)$ and observe that if there exists $k_0 \in \mathbb{N}$ such that $a_{k_0}=0$, then (ii) holds. Suppose now that there exists no $k_0 \in \mathbb{N}$ with $a_{k_0} = 0$, i.e., $a_k >0$ for all $k \in \mathbb{N}$. Then it follows from \eqref{comparacaok+1} that
\begin{equation}\label{lemmatecnical1}
a_{k+1} +\sigma a_{k+1}^p \leq a_k, \quad k\in \mathbb{N}.
\end{equation}
Remembering that  $a_k>0$ for all $k\in \mathbb{N}$, we divide both sides of \eqref{lemmatecnical1} by $a_{k+1}$  and take the limit as $k\longrightarrow \infty$ therein. This gives us
$$
\lim_{k\longrightarrow\infty}\left( 1+ \sigma a_{k+1}^{p-1}\right) \leq \lim_{k\longrightarrow\infty} \frac{a_k}{a_{k+1}}.
$$
By $a_k\downarrow 0$ and $p-1 <0$, the obtained inequality yields
$\displaystyle \lim_{k\longrightarrow \infty} \frac{a_k}{a_{k+1}}=\infty$, which indicates that $\{a_k\}$ converges superlinearly to zero and thus completes the proof .
\end{proof}\vspace*{-0.03in}

Now we are ready to establish the first result on convergence rates of any generic algorithm satisfying ($\mathcal{H}1$) and ($\mathcal{H}2$) under the exponent PLK conditions. The theorem below concerns the {\em value convergence} with the most significant improvement of the known results in the case of lower exponents.\vspace*{0.05in}

\begin{theorem}[\bf value convergence rates under exponent PLK and relative error conditions]\label{convergenciafinitaproximal} In the setting of Theorem~{\rm\ref{theo1}}, let $\{x_k\}$ be a generic sequence of iterates satisfying {\rm($\mathcal{H}1$)} and {\rm($\mathcal{H}2$)}, and let the exponent PLK property of $\phi$ hold at $\bar{x}$ with $\varphi(t)=Mt^{1-q}$ for some $M > 0$ and $q\in [0,1)$. The following convergence rates are guaranteed for the value sequence $\{\phi(x^k)\}:$\vspace*{-0.05in}
\begin{itemize}
\item[\bf(i)]  If $q=0$, then  the sequence $\{\phi(x^k)\}$ terminates at $\phi(\Bar x)$ in finite steps.

\item[\bf(ii)] If $q\in (0, \frac{1}{2})$, then either $\{\phi(x^k)\}$ terminates at $\phi(\Bar x)$ in a finite number of steps, or $\{\phi(x^k)\}$ converges superlinearly to $\phi(\Bar x)$ as $k\to\infty$. 

\item[\bf(iii)] If $q=1/2$, then $\{\phi(x^k)\}$ linearly converges to $\phi(\bar{x})$ as $k\to\infty$.

\item[\bf(iv)] If $\displaystyle q \in \left(1/2,1\right)$, then there exists a positive constant $\eta$ such that
\begin{equation}\label{iv}
\phi(x^k) - \phi(\Bar{x})  \leq \eta k^{-\frac{1}{2q-1}} \quad \text{for all large }\;k\in\mathbb{N}.
\end{equation}
\end{itemize}
\end{theorem}\vspace*{-0.1in}
\begin{proof}
It follows from Theorem~\ref{theo1} that $x^k\to\Bar x$ as $k\to\infty$ and that $\Bar  x$ is an $M$-stationary point of $\phi$. Defining  $\psi(x):= \phi(x) - \phi(\Bar{x})$, we see that the sequence $\{\psi(x^k)\} \subset [0,\infty)$ is decreasing. Using \eqref{desigualdadeklND} with $\displaystyle \varphi(t)= Mt^{1-q}$ as $q\in[0,1)$ and $M>0$, it follows that
\begin{equation}\label{klaplicada}
\begin{array}{ll}
&\| w^{k+1}\|\geq \text{ dist}\big(0, \partial\phi(x^{k+1})\big)  \geq \displaystyle\frac{1}{M(1-q)} \big( \phi(x^{k+1}) - \phi(\bar{x}) \big)^q\\
&\text{ with }\;w^{k+1} \in \partial\phi(x^{k+1}),\quad k\in\mathbb N. 
\end{array}
\end{equation}

To verify (i), suppose on the contrary that the sequence $\{x^k\}$ is infinitely generated and deduce from \eqref{klaplicada} with $q=0$ that $\| w^{k+1} \|\ge M^{-1}$ for all large $k$. This tells us together with \eqref{fran} in $({\cal H}2)$ that $(Mb)^{-1}\le\| x^{k+1} - x^k \|$, which contradicts the fact that $\{x^k\}$ converges to $\bar{x}$ and thus justifies the statement in (i).

For the proof of (ii), subtract $\phi(\bar{x})$ from both sides of \eqref{desigualdade1kNDH} and get
$$
\psi(x^{k+1})  \le \psi(x^k) - a \|x^{k+1} - x^k\|^2,
$$
which being combined with \eqref{fran} yields
$$
\psi(x^{k+1})  \leq \psi(x^k) - \frac{a}{b^2}\|w^{k+1} \|^2.
$$
Using the obtained inequality together with \eqref{klaplicada} and the definition of $\psi$ gives us
\begin{equation}\label{desigualdadegeral}
\psi(x^{k+1}) \leq \psi(x^k) - \frac{a}{b^2M^2(1-q)^2}\big(\psi(x^{k+1})\big)^{2q}.
\end{equation}
It follows from Lemma~\ref{lemmatec} with $\sigma =\frac{a}{b^2M^2(1-q)^2}$ and $p=2q$ that if $q\in (0,1/2)$, then the sequence $\{ \psi(x^k)\}$ is either finite or converges superlinearly, and thus (ii) is verified.

To prove (iii), let $q=1/2$ and $\alpha= \frac{4a}{b^2M^2}$. By \eqref{desigualdadegeral} we get
\begin{equation}\label{psi}
\psi(x^{k+1}) \leq \frac{1}{1+\alpha}\psi(x^k).
\end{equation}
Since $1/(1+\alpha) \in (0,1)$, assertion (iii) is a consequence of \eqref{psi}. 

It remains to verify (iv). For $q \in (1/2, 1)$, it follows from \eqref{desigualdadegeral} with $\beta:=\displaystyle\frac{a}{b^2M^2(1-q)^2}$ that we have
\begin{equation}\label{1}
\beta \leq  \frac{\psi(x^k) - \psi(x^{k+1}) }{\psi(x^{k+1})^{2q}}=[\psi(x^k) - \psi(x^{k+1})]\Phi(\psi(x^{k+1})),
\end{equation}
where $\displaystyle \Phi(t):= \frac{1}{t^{2q}}$. Suppose first that there exists  $Q\in (1,\infty)$ such that $\Phi(\psi(x^{k+1})) \leq Q \Phi(\psi(x^{k}))$. Combining this fact with \eqref{1} and Lemma~\ref{propositionAbsil2005} brings us to
\begin{eqnarray*}
\beta  &\leq & \frac{\psi(x^k) - \psi(x^{k+1}) }{\psi(x^{k+1})^{2q}}=\int_{\psi(x^{k+1})}^{\psi(x^k)}\Phi(\psi(x^{k+1}))dt \leq Q\int_{\psi(x^{k+1})}^{\psi(x^k)}\Phi(\psi(x^{k}))dt\, \nonumber \\  
& \leq & \frac{Q}{1-2q} \big[\big(\psi(x^k)\big)^{1-2q}- \big(\psi(x^{k+1})\big)^{1-2q}\big]= \frac{Q}{2q-1}\big[\big(\psi(x^{k+1})\big)^{1-2q} -\big(\psi(x^k)\big)^{1-2q}\big],
\end{eqnarray*}
which implies by $2q-1 > 0$ that
\begin{equation}\label{3}
\frac{\beta(2q-1)}{Q} \leq \big(\psi(x^{k+1})\big)^{1-2q} -\big(\psi(x^k)\big)^{1-2q}, \quad k\geq k_0. 
\end{equation}
Letting $j-1>k_0$, we deduce from \eqref{3} that
$$
\sum_{k=k_0}^{j} \frac{\beta(2q-1)}{Q} \leq 
\sum_{k=k_0}^{j} \big(\psi(x^{k+1})\big)^{1-2q} -\big(\psi(x^k)\big)^{1-2q}
$$
giving us in turn that
$$
\frac{(j-1-k_0)\beta(2q-1)}{Q} \leq \big(\psi(x^{j})\big)^{1-2q} - \big(\psi(x^{k_0})\big)^{1-2q}.
$$
Therefore, we arrive at the estimate
$$
\big(\psi(x^{j})\big)^{1-2q} \geq \frac{(j-1-k_0)\beta(2q-1)}{Q} + \big(\psi(x^{k_0})\big)^{1-2q}.
$$
Since the function $ t\mapsto t^{-\frac{1}{2q-1}}$ is decreasing for $q\in (1/2,1)$, it follows that
$$
\psi(x^{j}) \leq \left(\frac{(j-1-k_0)\beta(2q-1)}{Q} + \big(\psi(x^{k_0})\big)^{1-2q}\right)^{-\frac{1}{2q-1}},
$$
which ensures the existence of some $\eta >0$ with
$$
\psi(x^{j}) \leq \eta j^{-\frac{1}{2q-1}}.  
$$
and thus verifies the convergence rate estimate in \eqref{iv} in this case.

Consider now the alternative situation, where the number $Q$ used above does not exist, i.e., whenever $Q\in(1,\infty)$ we have $\Phi(\psi(x^{k+1}))> Q \Phi(\psi(x^{k}))$ for all $k$ sufficiently. Fix some $Q\in(1,\infty)$ and define $Q_1: = \frac{1}{Q^{\frac{1}{2q}}}$. Then it follows from the definitions that $\psi(x^{k+1})\leq Q_1 \psi(x^k)$ and hence $\psi(x^{k+1})^{1-2q}\geq Q_1^{1-2q}\psi(x^k)^{1-2q}$ by $1-2q<0$.  Subtracting $\psi(x^k)^{1-2q}$ from both sides of the latter inequality gives us
$$
\psi(x^{k+1})^{1-2q} - \psi(x^{k})^{1-2q} \geq \big(Q_1^{1-2q} -1\big)\psi(x^k)^{1-2q}. 
$$
Note that $Q_1 \in (0,1)$ and therefore $Q_1^{1-2q}>1$ for all $q\in(1/2,1)$. The convergence $\psi(x^k)\to 0$ yields  $\psi(x^k)^{1-2q}\to\infty$ as $k\to\infty$, and hence
$$
\psi(x^{k+1})^{1-2q} - \psi(x^{k})^{1-2q} \geq \frac{\beta(2q-1)}{Q}
$$
whenever $k$ is sufficiently large. The rest of the proof is similar to the above arguments based on  \eqref{3}, and thus we are done with verifying (iv) and the whole theorem.
\end{proof}\vspace*{-0.05in}

The next theorem establishes the convergence rates for the sequence of iterates $\{x^k\}$ of any algorithm satisfying ($\mathcal{H}1$) and ($\mathcal{H}2$) under the exponent PLK conditions. \vspace*{0.05in}

\begin{theorem}[\bf convergence rates of iterates under exponent PLK and relative error conditions]\label{convergenciafinitaproximalx^k} In the setting of Theorem~{\rm\ref{convergenciafinitaproximal}}, the following assertions hold:\vspace*{-0.05in}
\begin{itemize}
\item[\bf(i)] If $q=0$, then  the sequence $\{x^k\}$ terminates at $\Bar x$ in finite steps.

\item[\bf(ii)] If $q\in (0, \frac{1}{2})$, then either $\{x^k\}$ terminates at $\Bar x$ in finite steps, or 
$\{x^k\}$ converges superlinearly to $\Bar x$ as $k\to\infty$.

\item[\bf(iii)] If $q=1/2$, then $\{x^k\}$ converges linearly to $\bar{x}$ as $k\to\infty$.

\item[\bf(iv)] If $\displaystyle q \in \left(1/2,1\right)$, then there exists a positive constant $\sigma$ such that 
$$
\left\|x^k - \Bar{x} \right\| \leq \sigma k^{-\frac{1-q}{2q-1}}\;\text{ for all $k$ sufficiently large}.
$$
\end{itemize}
\end{theorem}\vspace*{-0.1in}
\begin{proof}
By taking into account that 
$$
0\leq a\|x^{k+1} - x^k \|^2 \leq \phi(x^k) - \phi(x^{k+1}), \quad k\in\mathbb{N},
$$ 
assertion (i) of this theorem is an immediate consequence of Theorem~\ref{convergenciafinitaproximal}(i).
To verify (ii), deduce first from the concavity of $\varphi$ in \eqref{desigualdadeklND} together with ($\mathcal{H}1$) and ($\mathcal{H}2$) that
\begin{equation}\label{11}
\frac{\|x^{k+1} - x^k\|^2}{\|x^k - x^{k-1}\|}\leq \frac{b}{a}\left[\varphi\left(\phi(x^k)-\phi(\bar x)\right) - \varphi\big(\phi(x^{k+1}) - \phi(\bar{x})\big) \right].
\end{equation}
for all large $k\in\mathbb N$. Given $r\in (0, 1)$, we have the two possibilities for such $k$:
\begin{itemize}
\item[\bf(a)] $\|x^{k+1} - x^k\| \geq r\|x^k - x^{k-1}\|$.
\item[\bf(b)] $\|x^{k+1} - x^k\| < r\|x^k - x^{k-1}\|$.
\end{itemize}
In case (a), the usage of \eqref{11} with $\varphi(t) = Mt^{1-q}$ and $\psi(x)=\phi(x)-\phi(\bar{x})$ yields
$$ 
\|x^{k+1} -x^k\| \leq \frac{bM}{ar}\big( \psi(x^k)^{1-q}  - \psi(x^{k+1})^{1-q}\big)\leq r \|x^k - x^{k-1} \|  +
\frac{bM}{ar}\big( \psi(x^k)^{1-q}  - \psi(x^{k+1})^{1-q}\big).
$$
On the other hand, in case (b) we readily get
$$
\| x^{k+1} -x^k \|\leq r \|x^k - x^{k-1} \| \leq r \|x^k - x^{k-1} \| + \frac{bM}{ar}\big( \psi(x^k)^{1-q}  - \psi(x^{k+1})^{1-q}\big),
$$
and therefore arrive at the inequality
\begin{equation*}
\| x^{k+1} -x^k \| \leq r \|x^k - x^{k-1} \| + \frac{bM}{ar}\big(\psi(x^k)^{1-q}-\psi(x^{k+1})^{1-q}\big).
\end{equation*}
The latter provides the following estimates valid for all the natural numbers $l\ge 1$:
\begin{equation}\label{somatorio}
\sum_{j=k}^{k+l}\| x^{j+1} -x^j \| \leq \frac{r}{(1-r)}\|x^k - x^{k-1} \| + \frac{bM}{ar(1-r)}\big( \psi(x^k)^{1-q}  - \psi(x^{k+l+1})^{1-q}  \big). 
\end{equation}
Define $s_k:= \sum_{j=k}^{\infty}\|x^{j+1} - x^j\|$ and let $l\to\infty$ in \eqref{somatorio}. Remembering that
$\{\psi(x^k)\}$ decreases to zero, we get from \eqref{somatorio} that
$$
s_k \leq \frac{r}{(1-r)}\|x^k - x^{k-1} \| + \frac{bM}{ar(1-r)}\psi(x^k)^{1-q}.
$$
Combining the last inequality with \textbf{($\mathcal{H}1$}) and using $\phi(x^{k-1})-\phi(x^k)\le\phi(x^{k-1})-\phi(\bar x)=\psi(x^{k-1})$ bring us to the estimate
\begin{equation}\label{eq:1111}
s_k \leq \frac{r}{a^{1/2}(1-r)} \big(\psi(x^{k-1})\big)^{1/2} + \frac{bM}{ar(1-r)}\psi(x^k)^{1-q}. 
\end{equation}
When $q\in(0,1/2)$, it follows from the above that
\begin{equation}\label{sk}
s_k \leq \frac{r}{a^{1/2}(1-r)}\big( \psi(x^{k-1})  \big)^{1/2} + \frac{bM}{ar(1-r)}\psi(x^k)^{1/2}.
\end{equation}
Since the sequence $\{\psi(x^k)\}$ is monotonically decreasing, we deduce from \eqref{sk} that
$$
s_k \leq C \big( \psi(x^{k-1})  \big)^{1/2}\;\mbox{ with }\; C:=\frac{a^{1/2}r^2+ bM}{ar(1-r)},
$$
which allows us to verify (ii) by using Theorem~\ref{convergenciafinitaproximal} and the estimate $\| x^k - \bar{x}\|\leq s_k$. Assertions (iii) and (iv) of this theorem follow directly from Theorem~\ref{convergenciafinitaproximal}(iii,iv)), respectively. Indeed, in case (iii) we set \( q = 1/2 \) in \eqref{eq:1111}, while (iv) with \( q \in (1/2, 1) \) is deduced from the inequality \( \big( \psi(x^{k}) \big)^{1/2} \leq \big( \psi(x^{k}) \big)^{1-q} \) held for large \( k \).
\end{proof}\vspace*{-0.03in}

The following example presents a family of {\em univariate smooth functions} for which the classical proximal algorithm, satisfying $({\cal H}_1)$ and $({\cal H}_2)$, does {\em not terminate in finite steps} under the PLK conditions with lower exponents.\vspace*{0.05in} 

\begin{example}[\bf failure of finite termination for the proximal algorithm under PLK lower exponents]\label{example101}
{\rm Consider the parametric class of functions $\phi_\alpha(x):=|x|^{1+\alpha}$ with $0<\alpha<1$; cf.\ \cite[page~22]{polyak1987}, where such functions were considered from different prospectives. Each function $\phi_\alpha$ is continuously differentiable, achieves its global minimum at $\bar x=0$, is $C^1$, and satisfies the exponent  PLK condition with $q=\frac{\alpha}{1+\alpha}\in (0,1/2)$. The classical proximal algorithm reads for $\phi_\alpha$ as follows:
\begin{equation}\label{proximal}
x_{k+1}= {\rm argmin}_{x\in \mathbb{R}^n}\Big\{ |x|^{1+\alpha} + \frac{1}{2\lambda_k}| x - x_k|^2\Big\},\quad \lambda_k>0,\;k\in\mathbb{N}.
\end{equation}
If we suppose that $\{x_k\}$ is finite, then there exists $k_0\in\mathbb{N}$ such that $(1+\alpha){\rm sign}(x_{k_0})|x_{k_0}|^{\alpha}=0$ and consequently $x_{k_0}=0$. Applying the stationary rule to the function in \eqref{proximal} gives us the equation
\begin{equation*}
0=(1+\alpha){\rm sign}(x_{k_0})|x_{k_0}|^{\alpha} + x_{k_0} - x_{k_0-1},
\end{equation*}
which yields $x_{k_0-1}=x_{k_0}=0$. Recursively, we arrive at $x_0=0$ and thus conclude that for any starting point $x_0\ne 0$ the sequence of iterates in \eqref{proximal} must be infinite. Theorem~\ref{convergenciafinitaproximalx^k} tells us that $x_k\to 0$ as $k\to\infty$ superlinearly.} 
\end{example}\vspace*{0.05in}

The next example illustrates the alternative situation in Theorem~\ref{convergenciafinitaproximalx^k},
where the proximal algorithm {\em terminates in finite steps} in contrast to Example~\ref{example101}.\vspace*{0.05in} 

\begin{example}[\bf finite step termination of the proximal algorithm]
{\rm Let 
\begin{equation}\label{prox}
\phi(x) = \begin{cases} -x, & x < 0, \\ x^{3/2}, & x \geq 0. \end{cases}
\end{equation}
We see that the function $\phi$ in \eqref{prox} achieves its global minimum at $\bar x=0$, where it is nondifferentiable, and bit satisfies the exponent PLK property with $q=1/3\in (0,1/2)$ at $\bar x$. The proximal algorithm for \eqref{prox} is written as 
\begin{equation}\label{eq.10102}
x_{k+1} ={\rm argmin}_{x \in \mathbb{R}^n}\Big\{ \phi(x) + \frac{1}{2\lambda_{k}} | x - x_{k} |^2\Big\}, \quad \lambda_k>0,\;k\in\mathbb{N}, 
\end{equation}
First we claim that if $x_0<0$, then $x_k \leq 0$ for all $k\in \mathbb{N}$. On the contrary, suppose that
there exists $k_0 \in \mathbb{N}$ for which $x_{k_0} < 0$ while $x_{k_0+1}>0$. By the smoothness of \eqref{prox} at the points under consideration, we deduce from the stationary rule applied to \eqref{eq.10102} that  
$$
x_{k_0} = \frac{3\lambda_{k_0}}{2} x_{k_0+1}^{\frac{1}{2}} + x_{k_0+1},
$$
which is a contradiction that verifies the claim. Letting  $\lambda_k := \bar{\lambda}>0$ for all $k\in \mathbb{N}$ and assuming  that the sequence is infinitely generated, we get $x_k < 0$. Applying again the stationary rule to  \eqref{eq.10102} in this setting yields
$\bar{\lambda} = x_{k+1} - x_k$.  which leads us to a contradiction. Therefore, for any fixed $\bar{\lambda}>0$, the sequence $\{x_k\}$ converges in finite steps if the starting point is selected as $x_0<0$.}
\end{example}\vspace*{0.05in}

We finalize this section with the following comments.\vspace*{0.05in}

\begin{remark}[\bf comparison with known results]\label{rem2} $\,$\vspace*{-0.05in} {\rm 
\begin{itemize}
\item[\bf(i)]
The main novelty of Theorem~\ref{convergenciafinitaproximal} and Theorem~\ref{convergenciafinitaproximalx^k} is for the case of lower exponents $q\in(0,1/2)$, where the obtained results significantly improve the known ones with {\em linear convergence} for various algorithms satisfying (${\cal H}_1$) and  (${\cal H}_2$); see \cite[Theorem~3.4]{frankel2015splitting} and many other publications in this direction.\footnote{When the results of our paper have been completed and prepared for publication, we got familiar with preprint \cite{pan} containing an alternative proof of superlinear convergence of descent algorithms satisfying $({\cal H}_1$) and (${\cal H}_2$) under the PLK conditions with lower exponents.}

\item[\bf(ii)] The convergence rate results obtained in Theorems~\ref{convergenciafinitaproximal}, \ref{convergenciafinitaproximalx^k} for the cases where $q=0$ and $q\in(1/2,1)$ under the exponent PLK conditions are closely related to those established in \cite[Theorem~3.4]{frankel2015splitting} for an implicit abstract model, which encompasses our explicit model defined by ($\mathcal{H}1)$ and ($\mathcal{H}2)$. The proofs of Theorems~\ref{convergenciafinitaproximal}, \ref{convergenciafinitaproximalx^k} are different from those given in \cite{frankel2015splitting}. 
\end{itemize}}
\end{remark}
\vspace*{-0.2in}


\section{Convergence Rates under Modified Error Property}\label{sec:rel-modif}\vspace*{-0.05in}


In this section, we provide a comprehensive rate convergence analysis of the generic algorithms satisfying (${\cal H}_1)$ and (${\cal H}_3)$ under the exponent PLK conditions. It is striking to observe that, in contrast to the results of Section~\ref{sec:rel-error}, the replacement of (${\cal H}_2)$ by (${\cal H}_3)$ fully {\em excludes 
infinite sequences of iterates} in the case of PLK with lower exponents. Moreover, we show below that, under some additional assumption, the lower exponent PLK property is {\em inconsistent} with the conditions in (${\cal H}_1)$ and (${\cal H}_3)$.\vspace*{0.05in}

\begin{theorem}[\bf convergence rates under exponent PLK and modified relative error conditions]\label{convergenciafinita} In the setting of Theorem~{\rm\ref{theo1}}, assume that the exponent PLK property of $\phi$ holds at $\bar{x}$ with $\varphi(t)=Mt^{1-q}$ for some $M > 0$ and $q\in [0,1)$. The following convergence rates are guaranteed for the generic iterative sequences satisfying the conditions in ${\rm({\cal H}_1)}$ and ${\rm({\cal H}_3)}:$
\begin{itemize}
\item[\bf(i)] If $q\in [0, \frac{1}{2})$, then $\{x^k\}$ and $\{\phi(x^k)\}$ are terminated in a finite number of steps to $\bar x$ and $\bar{\phi}=\phi(\bar x)$, respectively.

\item[\bf(ii)] If $q=\frac{1}{2}$, then the sequences $\{x^k\}$ and $\{\phi(x^k)\}$ converge linearly at $\Bar{x}$ and $\bar{\phi}$, respectively.

\item[\bf(iii)] If $\displaystyle q \in(1/2,1)$, then there exists a positive constant $\sigma$ such that
$$
\left\|x^k - \Bar{x} \right\| \leq \sigma k^{-\frac{1-q}{2q-1}} \quad \text{for all large }\;k\in\mathbb{N}.
$$
\end{itemize}
\end{theorem}\vspace*{-0.05in}
\begin{proof} It is shown in the proof of Theorem~\ref{theo1} that $\phi(x^k)\to\phi(\bar{x})$ as $k\to\infty$. Defining  $\psi(x):= \phi(x) - \phi(\Bar{x})$, we see that the sequence $\{\psi(x^k)\} \subset [0,\infty)$ is decreasing, and hence $\psi(x^j)\ge t$ for all $t\in [\psi(x^{j+1}), \psi(x^j)]$, $j\in\mathbb{N}$. Thus
\begin{equation*}
\frac{1}{\psi(x^j)^q}\le\frac{1}{t^q}\;\mbox{ whenever }\;t\in [\psi(x^{j+1}),\psi(x^j)]\;\mbox{ and }\; q \in (0,1),
\end{equation*}
which implies in turn that
\begin{equation*}
\frac{\psi(x^j) - \psi(x^{j+1})}{\psi(x^j)^q} =\int_{\psi(x^{j+1})}^{\psi(x^j)} \frac{1}{\psi(x^j)^q}dt\le\int_{\psi(x^{j+1})}^{\psi(x^j)}\frac{1}{t^q}dt.
\end{equation*}
Therefore, we arrive at the inequality
\begin{equation*}
\frac{\psi(x^j) - \psi(x^{j+1})}{\psi(x^j)^q}\le\frac{1}{1-q }\big(\psi(x^j)^{1- q} - \psi(x^{j+1})^{1-q}\big),
\end{equation*}
which gives us the estimate
\begin{equation}\label{relacao2ND}
\sum_{j=k}^{k+l} \frac{\psi(x^j) - \psi(x^{j+1})}{\psi(x^j)^q} \leq \frac{1}{1-q}\left(\psi(x^k)^{1-q} - \psi(x^{k+l+1})^{1-q}\right)\;\mbox{ for all }\;k,l \in \mathbb{N}.
\end{equation}
On the other hand, combining \eqref{desigualdade1kNDH}, the construction of $\psi$, and \eqref{fran2} in $({\cal H}3)$ yields
\begin{equation*}
\psi(x^{k+1})\le \psi(x^k) - a\|x^{k+1} - x^k\|^2\le \psi(x^k)-\frac{a}{b^2}\| w^k \|^2.
\end{equation*}
By the PLK property \eqref{desigualdadeklND} with  $\varphi(t)=Mt^{1-q}$, we find $M>0$ such that
\begin{equation}\label{DCArelacao1ND}
\| w^k \|\geq \text{dist}\big(0, \partial \phi(x^k)\big) \geq \frac{1}{M(1-q)} \big(\phi(x^k) - \phi(\Bar{x})\big)^q.
\end{equation}
Using the latter together with \eqref{DCArelacao1ND} tells us that
\begin{equation}\label{termonegativoND}
\psi(x^{k+1}) \leq \psi(x^k) - \frac{a}{b^2M^2(1-q)^2}\big(\psi(x^k)\big)^{2q}.
\end{equation}

To furnish now the proof of assertion (i), consider first the case where $q\in (0,1/2)$ and suppose, arguing by contradiction, that the sequence $\{x^k\}$ is infinitely generated, which tells us that $\psi(x^k)>0$ for all $k\in\mathbb{N}$. Recall the convergence $\psi(x^k)\to 0$ as $k\to\infty$ and deduce from $2q-1 < 0$, $\eqref{termonegativoND}$, and the equality
$$
\psi(x^k)- \frac{a}{b^2M^2(1-q)^2}\psi(x^k)^{2q} = \psi(x^k)\Big(1- \frac{a}{b^2M^2(1-q)^2}\Big)\big(\psi(x^k)\big)^{2q-1}
$$
that there exists a number $k_0 \in \mathbb{N}$ such that
$$
\psi(x^{k_0 +1}) \leq \psi(x^{k_0})\Big(1- \frac{a}{b^2M^2(1-q)^2}\big(\psi(x^{k_0})\big)^{2q-1}\Big)<0,
$$
which contradicts the fact that $\psi(x^k)>0$ for all $k\in \mathbb{N}$. Thus for $q\in(0,1/2)$, the sequence $\{x^k\}$ reaches the $M$-stationary point $\Bar{x}$ in a finite number of steps.

To verify (i), it remains to consider the case where $q=0$. Arguing again by contradiction, suppose that the sequence $\{x^k\}$ is infinitely generated and deduce from \eqref{DCArelacao1ND} with $q=0$ that $\| w^k \|\ge M^{-1}$ for all large $k$, which tells us together with \eqref{fran2} in $({\cal H}_3)$ that $(Mb)^{-1}\le\| x^{k+1} - x^k \|$. This contradicts the fact that $\{x^k\}$ converges to $\bar{x}$ and thus concludes the proof of assertion (i).

To verify assertion (ii), we get from \eqref{termonegativoND} with $q=1/2$ that
\begin{equation}\label{eq:1001ND}
\psi(x^{k+1}) \leq\Big( 1-\frac{4a}{b^2M}\Big)\psi(x^k).
\end{equation}
Due to $\psi(x^k)\geq 0$ as $k\in \mathbb{N}$, it follows from \eqref{eq:1001ND} that $4a(b^2M)^{-1}\in (0,1)$, and hence $\{\phi(x^k)\}$ converges linearly to
$\phi(\bar{x})$. On the other hand, we have by \eqref{desigualdade1kNDH} and \eqref{fran2} that
\begin{equation*}
\phi(x^{k+1})\leq \phi(x^k) - a b^{-1}\|x^{k+1} -x^k\|\cdot\| w^k\| \;\mbox{ for all }\;k\in\mathbb N,
\end{equation*}
which being combined with \eqref{DCArelacao1ND} and the construction of $\psi$ yields
\begin{equation*}
\psi(x^{k+1}) \leq \psi(x^k) - \frac{ 2a }{Mb}\| x^{k+1} - x^k \|\big(\psi(x^k)\big)^{\frac{1}{2}}.
\end{equation*}
The obtained inequality readily ensures the estimate
\begin{equation*}
\|x^{k+1} -x^k\|\le\frac{Mb}{2a}\Big(\frac{ \psi(x^k) - \psi(x^{k+1})}{\psi(x^k)^{\frac{1}{2}}}\Big),
\end{equation*}
which tells us together with \eqref{relacao2ND} that
\begin{eqnarray*}
\sum_{j=k}^{k+l} \|x^{j+1} - x^j\| &\leq & \frac{Mb}{2a} \sum_{j=k}^{k+l} \frac{ \psi(x^j) - \psi(x^{j+1})}{\psi(x^j)^{\frac{1}{2}}}\, \\ & \leq & \frac{Mb}{a}\big(\psi(x^{k})^{\frac{1}{2}} - \psi(x^{k +l+1})^{\frac{1}{2}}\big).
\end{eqnarray*}
By taking the limit in the last inequality as $l\to\infty$, we arrive at
\begin{equation*}
s_k:= \sum_{j=k}^{\infty} \|x^{j+1} - x^j\|\le\frac{Mb}{a}\big(\psi(x^k)\big)^{\frac{1}{2}},
\end{equation*}
which being combined with \eqref{DCArelacao1ND} and \eqref{fran2} brings us to the estimates
\begin{equation*}
s_{k+1} \leq \frac{Mb}{a}\big(\psi(x^{k})\big)^{\frac{1}{2}}\le\frac{M^2b}{2a}\| w^k \| \leq\frac{M^2b^2}{2a}\| x^{k+1} - x^k \|.
\end{equation*}
Since $\|x^{k+1} -x^k \| = s_k - s_{k+1}$, the last inequality implies that
$$
s_{k+1} \le(2a)^{-1}M^2b^2(s_k -s_{k+1})
$$
from which it clearly follows that
\begin{equation*}
s_{k+1} \leq \frac{\tau}{1+\tau} s_k\;\mbox{ with }\;\tau:=\frac{M^2b^2}{2a}.
\end{equation*}
This tells us that $\{s_k\}$ converges linearly to $0$ as $k\to\infty$. On the other hand, using the triangle inequality yields
$\|x^k - x^{k+l}\| \leq \sum_{j=k}^{k+l} \|x^{j+1} - x^j\|$
from which we get that $\|x^k - \bar{x} \|\leq s_k$. Therefore, the linear convergence of $x^k\to\bar{x}$ follows from the linear convergence of $s_k\to 0$, and thus (ii) is fully justified. Assertion (iii) can be verified similarly to the proof of \cite[Theorem~2]{attouch2009convergence} given in the case of the proximal algorithm.
\end{proof}

The next significant result shows that the lower exponent PLK property of the cost function $\phi$ is {\em inconsistent} {\em fails} with $({\cal H}1)$ and $({\cal H}3)$ under the fulfillment of some additional assumption on the generic algorithm iterates. Note that another inconsistency result is established in Section~\ref{sec:KL} for a class of difference programming problems under the assumptions independent of $({\cal H}1)$ and $({\cal H}3)$.\vspace*{0.05in}

\begin{theorem}[\bf inconsistency of lower exponent PLK with  $(\bf{\cal H}1)$ and $(\bf{\cal H}3)$]\label{PLKfail} In addition to the general assumption of Theorem~{\rm\ref{convergenciafinita}}, suppose that 
$x^{k+1} \neq x^k$ for all $k\in \mathbb{N}$. Then the exponent PLK property of $f$ at $\bar x$ fails whenever $q\in(0,1/2)$.
\end{theorem}\vspace*{-0.05in}
\begin{proof}
Observe that if $x^k \neq x^{k+1}$ whenever $k\in \mathbb{N}$, then it follows from
$$
a \|x^{k+1} - x^k\|^2 \le\phi(x^k) - \phi(x^{k+1})
$$
that $\phi(x^k) \neq \phi(x^{k+1})$ for all $k \in \mathbb{N}$. Arguing by contradiction, suppose that $\phi$ satisfies the PLK property at $\bar{x}$ with exponent $q\in (0, 1/2)$. It follows from $({\cal H}1)$ and $({\cal H}3)$ that
$$
\phi(x^k)-\phi(x^{k+1}) \geq \frac{a}{b^2}\|w^k\|^2.
$$
Combining the last inequality with the the exponent PLK condition yields
$$ 
\phi(x^k)-\phi(x^{k+1}) \geq \frac{a}{(1-q)^2M^2b^2}\big(\phi(x^k)-\phi(\bar x) \big)^{2q}.
$$
Using the latter together with $\phi(x^k)-\phi(\bar x) \geq \phi(x^k)-\phi(x^{k+1})$ brings us to
\begin{equation}\label{beta}
\left(\phi(x^k) - \phi(\bar{x})\right)^{1-2q} \geq \beta,\;\mbox{ where }\;\beta:= \frac{a}{(1-q)^2M^2b^2}.
\end{equation}
Since $\phi(x^k) - \phi(\bar{x})$ as $k\to\infty$, \eqref{beta} fails for large $k$ when
$q \in (0,1/2)$. The obtained contradiction verifies the claimed assertion and completes the proof of theorem.
\end{proof}\vspace*{-0.03in}

\begin{remark}[\bf novelty under the modified elative error condition] $\;$ \rm
\begin{itemize}\vspace*{-0.05in}
\item[\bf(i)] Similarly to the case of the relative error condition (${\cal H}_2$) discussed in Remark~\ref{rem2}, the convergence rate results obtained in Theorem~\ref{convergenciafinita} under the modified version ({\cal H}3) and the exponent PLK properties with $q=0$ and $q\in[1/2,1)$ continue to be generic extensions of those given for particular algorithms (see, e.g., \cite{Aragon2018,aragon2023coderivative,aragon2020boosted,attouch2009convergence,frankel2015splitting,li,Mordukhovich2023inexact,mordukhovich2024}).

\item[\bf(ii)] 
In the case of {\em continuous-time} gradient/subgradient systems, a {\em finite-time} convergence estimate has already appeared in the literature when the exponent PLK varies within the interval (0,1/2); see \cite[Theorem~2.7]{chill2006convergence}. Inspired by the latter result, it is shown in  \cite[Theorem~3.5(i)]{frankel2015splitting} that for certain algorithms including gradient-related methods, this better estimate remains true at least when dealing with the sequence of functional values. Our approach to this issue is novel being also in accordance with the estimates presented in the continuous case.

\item[\bf(iii)] To the best of our knowledge, the {\em inconsistency} result in the setting of Theorem~\ref{PLKfail} has never been observed before. Its analogs in the different setting are discussed in Section~\ref{sec:KL}; see Theorem~\ref{lipDC} and the discussion after its proof.
\end{itemize}
\end{remark}\vspace*{-0.2in}

\section{Specifications of Generic Descent Algorithms}\label{sec:specif}\vspace*{-0.05in}

This section discusses some specifications of the general class of descent methods and convergence analysis for them developed in Sections~\ref{sec:conv} and Section~\ref{sec:rel-modif}.\vspace*{0.03in} 

We consider here two particular methods and start with the following model.\vspace*{-0.1in}

\subsection{Inexact Reduced Graduate Methods}\vspace*{-0.05in}

A class of {\em inexact reduced graduate $($IRG$)$ methods} with various stepsize selections have been recently proposed and developed by Khanh et al. \cite{Mordukhovich2023inexact} for problems of smooth nonconvex optimization. The model unifying the IRG methods in \cite{Mordukhovich2023inexact} considers the objective function $\phi\colon\mathbb{R}^n\to\mathbb R$ and the iterative sequence
$\{x^k\}$ as $k\in\mathbb N$ satisfying
\begin{equation}\label{relacaoabstrata}
\phi(x^k) -\phi(x^{k+1}) \geq \frac{\beta}{t_k} \|x^{k+1} - x^k \|^2\;\text{ and } \;\|\nabla\phi(x^k)\| \leq \frac{c}{t_k}\|x^{k+1} -x^k\|,\quad
\end{equation}
where ${t_k}\subset \mathbb{R}_+$ and $\beta, c>0$. This can be viewed as a smooth version of the general nonsmooth scheme investigated above. The convergence analysis of IRG algorithms conducted in \cite{Mordukhovich2023inexact} is based on the observation that the sequences of iterates generated by \eqref{relacaoabstrata} satisfy the following properties:\vspace*{-0.05in}
\begin{itemize}
\item \textit{Primary descent condition}: There exists $\sigma > 0$ such that 
$$
\phi(x^k)-\phi(x^{k+1}) \ge\sigma \|\nabla\phi(x^k) \|\cdot\| x^{k+1}- x^k\|\;\mbox{ for large }\;k\in\mathbb N.
$$
\item \textit{Complementary descent condition}: we have the implication
$$
\big[\phi(x^{k+1})=\phi(x^k)\big] \Longrightarrow\big[x^{k+1}=x^k\big]\;\mbox{ for large }\;k\in\mathbb N.
$$
\end{itemize}
In fact, the latter conditions, while not \eqref{relacaoabstrata}, are used by Absil et al. in \cite[page~536]{Absil2005} to prove the convergence of the corresponding iterative sequences when the cost function is real analytical. Let us emphasize that the general scheme \eqref{relacaoabstrata} comprises {\em new linesearch methods} with an {\em inexact gradient information} for finding stationary points of ${\cal C}^1$-smooth functions by using different choices of stepsize rules. To the best of our knowledge, \cite{Mordukhovich2023inexact} is a first paper dealing with convergence rates for methods satisfying condition \eqref{relacaoabstrata}. In particular, it extends some special settings considered, e.g., in \cite{polyak1963gradient} and \cite{karimi2016linear}, where the authors analyze convergence properties of the exact gradient descent method for functions satisfying the PLK conditions. More precisely, \cite[Theorem~2.5]{Mordukhovich2023inexact} extends the scope of applicability of the results in the aforementioned papers. The results established above allow us to essentially improve those obtained in  \cite{Mordukhovich2023inexact} for IRG methods and their various linesearch specifications. In particular, we now have the {\em finite termination} of such algorithms under \eqref{relacaoabstrata} and  the exponent PLK with $q\in[0,1/2)$ at any {\em accumulation point} of the iterative sequences {\em provided} that the corresponding stepsize sequence $\{t_k\}$ is {\em bounded away} from zero. The latter assumption is required to accommodate $(\rm{\cal H}3)$ with a fixed $b>0$ in \eqref{fran2}. \vspace*{-0.1in}

\subsection{Boosted DC Algorithm}\label{boosted}\vspace*{-0.05in}

Let us recall the following Boosted DC Algorithm with Backtracking (BDCA), which was proposed and investigated by Arag\'on-Artacho and Vuong \cite{aragon2020boosted} to solve the problem:
$$
(\mathcal{P}) \qquad \min_{x\in\mathbb{R}^n} \phi(x):= g(x) -h(x),
$$
where both $g,h:\mathbb{R}^n \to\overline{\mathbb{R}}$ are  assumed to be strongly convex with modulus $\rho>0$, with $g$ being strictly differentiable at the points in question:\vspace*{-0.15in}

\begin{algorithm}[H] \caption{BDCA}
\begin{algorithmic}[1]
\item[1.] Fix $\alpha>0, \Bar{\lambda}>0$ and $ \beta\in (0,1)$. Let $x_0$ be any initial point and set $k:= 0$.
\item[2.]  Select $u_k \in \partial h(x^k)$ and solve the strongly convex optimization problem
$$
(\mathcal{P}_k')\,\,\,\,\,\,\,\,\,\min_{x\in \mathbb{R}^n} g(x) - \langle u_k,, x \rangle
$$
to obtain its unique solution $y^k$.
\item[3.] Let $d^k:=  y^k - x^k$. If $d^k = 0$, \textbf{stop} and \textbf{return} $x^k$. Otherwise, go to Step~4.
\item[4.] Choose any $\Bar{\lambda}_k\geq 0$. Set $\lambda_k:= \Bar{\lambda}_k$.
\textbf{ While } $\phi(y^k +\lambda_kd^k) > \phi (y^k) - \alpha \lambda_{k}^2\|d^k\|^2$, \textbf{ do } $\lambda_k:= \beta \lambda_k.$
\item[5.] Let $x^{k+1} := y^k + \lambda_k d^k$. If $x^{k+1} = x^k$, \textbf{ stop} and \textbf{return} $x^k$. Otherwise, set $ k := k+1 $ and go to Step~2.
\end{algorithmic}
\end{algorithm}\vspace*{-0.15in}
When $h$ is also smooth, BDCA was introduced and analyzed in Arag\'on-Artacho et al. \cite{Aragon2018}. In both papers, BDCA accelerates the convergence of the classical Difference of Convex Functions Algorithm (DCA). It is shown in \cite{aragon2020boosted} that the {\em strong PLK condition} imposed on $\phi$ at an accumulation point yields the global convergence of iterates with deriving convergence rates under the exponent strong PLK condition with $\varphi(t)= Mt^{1-q}$ for $q\in[0,1)$ and $M>0$. Moreover, it is observed in \cite[Remark~4.4]{aragon2020boosted}, with the reference to the second author of the current paper, that for Lipschitz continuous functions $\phi$, the strong PLK property may be replaced by the symmetric one.

Next we show that BDCA can be viewed as a specification, for non-Lipschitzian {\em continuous} functions $\phi$, of the generic algorithm developed in Section~\ref{sec:conv} with using the {\em symmetric subdifferential} \eqref{sym} and the {\em symmetric PLK property} of $\phi$. This approach allows us to improve and extend the results of \cite{aragon2020boosted} about the convergence and convergence rates of BDCA to the case of continuous functions. First we recall some properties of the symmetric subdifferential of continuous functions used below; all of them and much more on \eqref{sym} (including full calculus) can be found in \cite{mordukhovich2006variational,mordukhovich2018variational}.\\[1ex]
\noindent $\bf(\mathcal{A}_1)$  $\partial^0(-\phi(x)) = -\partial^0 \phi(x)$ (plus-minus symmetry).\\[1ex]
\noindent $\bf(\mathcal{A}_2)$ $\partial^0 \phi(x) = \{ \nabla \phi (x) \}$ when $\phi$ is strictly differentiable.\\[1ex]
\noindent $\bf(\mathcal{A}_3)$ We have the representation $\partial^0 \phi(\bar x) = \partial \phi(\bar x) = \big\{x^* \in \mathbb{R}^n\;|\;\langle x^*, x - \bar{x} \rangle \leq \phi(x) - \phi(\bar{x})\;\mbox{ for all }\; x \in \mathbb{R}^n\}$ when $\phi$ is convex.\\[1ex]
\noindent $\bf(\mathcal{A}_4)$ If $\phi=f_1+f_2$ with $f_1$ being strictly differentiable, then  we have the sum rule $\partial^0\phi(x) = \nabla f_1(x) + \partial^0 f_2(x)$.\vspace*{0.05in}

Observe that $\partial^0\psi(x)$ is often nonconvex being significantly smaller than Clarke's subdifferential $\overline{\partial}\psi(x)$, with the relationship $\overline{\partial}\psi(x)={\rm co}\,\partial^0\phi(x)$ when $\phi$ is locally Lipschitzian around $x$. Note also that Clarke's counterpart of the plus-minus symmetry property as in $(\mathcal{A}_1)$  requires the local Lipschitz continuity of $\phi$ around $\bar x$.\vspace*{0.03in}

The theorem below incorporates BDCA into the generic class of algorithms satisfying conditions (${\cal H}1)$ and  (${\cal H}3)$ via the symmetric subdifferential and provides its convergence analysis under the symmetric exponent PLK properties.\vspace*{0.05in}

\begin{theorem}[\bf BDCA for continuous functions] Consider problem $({\cal P})$, where $\phi$ is continuous on its domain. Let the sequence $\{x^k\}$ be generated by BDCA, and let $\nabla g$ be Lipschitz continuous with modulus $L$ around an accumulation point $\bar x$ of
$\{x^k\}$. Then for all large $k\in\mathbb{N}$ we have the conditions:
\begin{equation}\label{desigualdade1NDdiscution}
\phi(x^{k+1})\leq \phi(x^k) - \frac{\alpha\lambda_k^2 +\rho}{(1+\lambda_k)^2} \|x^{k+1} -x^k\|^2,
\end{equation}
\begin{equation}\label{desigualdade2NDdiscution}
\mbox{there exists }\;w^k \in \partial^0 \phi(x^k)\;\mbox{ such that }\;\| w^k\|\le L\| x^{k+1} -x^k\|,
\end{equation}\vspace*{0.03in}
and $\{x^k\}$ converges to $\bar x$ as $k\to\infty$, which is a symmetric stationary point of $({\cal P})$ satisfying $0\in\partial^0\phi(\bar x)=\nabla g(\bar x)-\partial h(\bar x)$. If furthermore $\phi$ enjoys the exponent PLK property with $\varphi(t)=M t^{1-q}$ for some $M>0$ and $q\in[0,1)$, then the convergence rates of $x^k\to\bar x$ are as given in Theorem~{\rm\ref{convergenciafinita}}.
\end{theorem}\vspace*{-0.05in}
\begin{proof} Using the definitions of $x^{k+1}$ and $d^k$ in BDCA gives us the equalities
\begin{equation*}
\| x^{k+1} -x^k\|=\|y^k +\lambda_k d^k -x^k\|= \| y^k+ \lambda_k(y^k-x^k) -x^k \| =(1+\lambda_k)\|d^k\|
\end{equation*}
from which we deduce the estimate
\begin{equation}\label{preliminar1}
\|d^k\|=\frac{1}{1+\lambda_k}\| x^{k+1} -x^k\|\leq \| x^{k+1} -x^k\|.
\end{equation}
Hence the estimate in \eqref{desigualdade1NDdiscution} follows from the conditions
$$
\phi (y^k)\leq \phi(x^k) -\rho \|d^k\|^2\;\mbox{ and }\;
\phi(y^k +\lambda_k d^k)\leq \phi(y^k) - \alpha \lambda^2_k\|d^k\|^2\;\mbox{ for some }\;\delta_k>0
$$
obtained for BDCA in \cite[Proposition~3.1]{aragon2020boosted}. To verify further \eqref{desigualdade2NDdiscution}, we employ properties $({\cal A}_1)-({\cal A}_4)$ of the symmetric subdifferential and get the relationships
$$
\nabla g(y_k) - \nabla g(x_k) \in \partial h(x_k) - \nabla g(x_k) = \partial^0(-\phi(x_k)) = -\partial^0\phi(x_k).
$$
Recalling that  $y^k$ is the unique solution to the subproblem $(\mathcal{P}_k')$ in BDCA, i.e., $u^k= \nabla g(y^k)$ ensures that $w^k := \nabla g(x^k) - \nabla g(y^k)\in \partial^0\phi(x^k)$. Take now a neighborhood $U$ of $\bar x$ on which $\nabla g$ is Lipschitz continuous with modulus $L$. It follows from \eqref{preliminar1} that $x^k, y^k\in U$ for all large $k\in\mathbb N$, and hence
\begin{equation*}
\|\nabla g(x^k)-\nabla g(y^k)\|\le L\|x^k-y^k\|,
\end{equation*}
which yields \eqref{desigualdade2NDdiscution}. Since $\bar{x}$ is an accumulation point of $\{x^k\}$ and $\phi$ is continuous on its domain, we deduce from \eqref{desigualdade1NDdiscution} and \eqref{desigualdade2NDdiscution} by the proof of Theorem~\ref{theo1} that $x_k\to\bar x$ as $k\to\infty$ with $0\in\partial^0\phi(\bar x)$. The latter implies by $({\cal A}_1)-({\cal A}_4)$ that $0\in\nabla g(\bar x)-\partial h(\bar x)$. The convergence rates under the exponent PLK conditions follows from Theorem~\ref{convergenciafinita}.
\end{proof}\vspace*{-0.05in}

Note that the finite termination of BDCA at any accumulation point under the exponent PLK condition with $q\in(0,1/2)$ is a significant improvement of \cite{aragon2020boosted} even for Lipschitz continuous objectives, where merely linear convergence of $\{x^k\}$ is justified. We'll see in Section~\ref{sec:KL} the further clarification of this fact at local minimizers of $\phi$.\vspace*{-0.1in}

\section{Lower Exponent PLK for Difference Programs}\label{sec:KL}\vspace*{-0.05in}

In this section, we take a close look at the exponent PLK property from Definition~\ref{def:kl}(iv) with {\em lower exponents} $q\in(0,1/2)$ for a class of {\em difference programs}.

As shown in \cite[Theorem~1]{fee}\footnote{The authors are grateful to Guoyin Li for drawing our attention to paper \cite{fee}.}, the classical PL gradient inequality always holds for any {\em analytic} function $\phi\colon\mathbb{R}^n\to\mathbb R$ with some exponent $q\in[1/2,1)$. Although it does not exclude the fulfillment of the PLK conditions with lower exponents, we show below that this property {\em fails} for {\em local minimizers} of smooth functions with {\em Lipschitzian gradients}, i.e., of class ${\cal C}^{1,1}$, without any connections with decent algorithms satisfying (${\cal H}_1$) and (${\cal H}_3$) as in Theorem~\ref{PLKfail}. In fact, our analysis provides such a result for problems of {\em difference programming} of type $({\cal P})$ (see Subsection~\ref{boosted}) with $\phi=g-h$, where $g$ is of class ${\cal C}^{1,1}$ (may not be convex), while $h\colon\mathbb{R}^n\to\overline{\mathbb R}$ is l.s.c.\ and convex. Recall that for any $x\in{\rm int}({\rm dom}\,h)$, the directional derivative of $\phi$ is represented by
\begin{equation}\label{dirder}
\phi'(x, d) = \langle \nabla g(x), d\rangle - \max_{v \in \partial h(x)} \langle v, d \rangle\;\mbox{ whenever }\;d\in\mathbb{R}^n.
\end{equation}

Let us introduce the two important {\em measures of stationarity} at the point $x\in{\rm int}({\rm dom}\,h)$ as follows. Observe first from \eqref{dirder} that
\begin{eqnarray*}
\min_{\|d\|=1} \phi'(x, d) & = & \min_{\|d\|=1} \Big\{ \langle\nabla g(x), d \rangle -\max_{v\in\partial h(x)}\langle v, d\rangle \Big\} =\min_{\|d\|=1} \min_{v \in \partial h(x)} \langle\nabla g(x) - v, d\rangle \\
& = & \min_{v\in\partial h(x)}\big(\{ -\|\nabla g(x) - v\|\big\} = - \max_{v\in\partial h(x)} \| \nabla g(x) - v\| := -\sigma_1(x).
\end{eqnarray*}
Thus $\sigma_1(x)$ measures the {\em fastest local decrease} of the function $\phi$. For each $x\in \mathbb{R}^n$, fix
$$
v_1(x) \in {\rm argmax}_{v\in \partial h(x)} \| \nabla g(x) - v\|
$$
and note that this point satisfies the necessary condition for local maximum:
$$
\langle v_1(x) - \nabla g(x), v - v_1(x) \rangle \leq 0\;\mbox{ whenever }\; v \in \partial h(x).
$$
In other words, we have the equality
\begin{equation}\label{vazio2}
\langle \nabla g(x) - v_1(x), v_1(x) \rangle = \min_{v  \in \partial h(x)} \langle \nabla g(x) - v_1(x), v \rangle.
\end{equation}
The second measure of stationarity at $x$ ids defined from the conditions
\begin{eqnarray*}
\max_{\|d\|=1} \phi'(x,d) & = & \max_{\|d\|=1} \Big\{\langle\nabla g(x), d \rangle - \max_{v\in\partial h(x)} \langle v, d \rangle \Big\} = \max_{\|d\|\leq 1} \min_{v\in\partial h(x)} \langle \nabla g(x) -v, d\rangle \\
& = & \min_{v\in \partial h(x)} \max_{\|d\|\leq 1} \langle \nabla g(x) -v, d \rangle = \min_{v\in\partial h(x)} \| \nabla g(x) - v\| := \sigma_0(x).
\end{eqnarray*}
Thus $\sigma_0(x)$ measures the {\em fastest local increase} of $\phi$. We obviously get
\begin{equation*}
\sigma_0(x)\leq \sigma_1(x)\;\mbox{ whenever }\;x\in \mathbb{R}^n.
\end{equation*}

Given an arbitrary point $x_0 \in \mathbb{R}^n$, consider its {\em own level set} $\mathcal{L}(x_0)$, which is the maximal connected component in $\{x\;|\;\phi(x) \leq \phi(x_0)\}$ containing $x_0$. Let us define the {\em exact Lipschitzian bound} of $\nabla g$ with respect to ${\cal L}(x_0)$ by
\begin{equation}\label{exact}
L(x_0) = \inf\big\{L\;\big|\; \| \nabla g(x) - \nabla g(y) \| \leq L \| x - y \|\;\mbox{ over }\; x, y \in \mathcal{L}(x_0)\big\}.
\end{equation}

The next lemma measures the difference between $\phi(x_0)$ and the global minimum of $\phi$ with respect to the level set ${\cal L}(x_0)$  in terms of the exact Lipschitzian bound \eqref{exact} and the stationarity measure $\sigma_1(x_0)$.\vspace*{0.05in}

\begin{lemma}[\bf distance estimate from global minimum]\label{lemma2}
Let $\bar x\in{\rm int}({\rm dom}\,h)$ be a global minimizer of the function $\phi$ in $({\cal P})$ with respect to the level set $\mathcal{L}(x_0)$, i.e., 
$$
\phi(x) \geq \phi(\bar x)\;\mbox{ for all }\;x\in \mathcal{L}(x_0).
$$
If $L(x_0) <\infty$, then for $d(x_0): = \nabla g(x_0) - v_1(x_0)$ and any $\displaystyle\alpha \in [0, 1/L(x_0)]$, we have 
\begin{equation}\label{vazio4}
\phi(x_0) - \phi(\bar{x}) \geq \phi\big(x_0) - \phi(x_0 - \alpha d(x_0)\big) \geq \frac{1}{2}\alpha\big(2 - \alpha L(x_0)\big)\sigma_{1}^{2}(x_0).
\end{equation}
\end{lemma}\vspace*{-0.1in}
\begin{proof} Since the case where $\sigma_{1}(x_{0})=0$ is obvious, suppose that
$\sigma_{1}(x_{0}) > 0$. Denote $x_{\alpha}: = x_{0} - \alpha d(x_{0})$ for $\alpha \geq 0$ and define the function $\xi(\alpha): = \phi(x_{\alpha})$. Note that
\begin{eqnarray*}
\xi'(0, 1) & = &\phi'\big(x_{0}, -d(x_{0})\big) = \langle \nabla g(x_{0}), -d(x_{0})\rangle - \max_{v \in \partial h(x_{0})} \langle v, -d(x_{0})\rangle \\
& = & \langle \nabla g(x_{0}), -d(x_{0})\rangle + \min_{v \in \partial h(x_{0})} \langle v, d(x_{0})\rangle \stackrel{\eqref{vazio2}}{=} -\sigma_{1}(x_{0}) < 0.
\end{eqnarray*}
Let $\bar{\alpha}$ be the smallest element of the set $\Sigma = \{\alpha \geq 0\;|\; \xi'(\alpha, 1) \geq 0\}$. If $\Sigma = \emptyset$, then $\bar{\alpha} =\infty$. Since $\xi'(0, 1) < 0$, we have $\xi'(\alpha, 1) < 0$ for all $\alpha \in [0, \bar{\alpha})$. Thus any $x_{\alpha}$ with $\alpha \in [0, \bar{\alpha}]$ belongs to the set $\mathcal{L}(x_{0})$. Adding the fact that $\partial h$ is monotone, this yields
\begin{equation*}
\begin{aligned}
\xi'(\alpha, 1) & = \langle \nabla g(x_{\alpha}), -d(x_{0})\rangle - \max_{v\in \partial h(x_{\alpha})} \langle v, -d(x_{0})\rangle \\
& \leq \langle \nabla g(x_{\alpha}), -d(x_{0})\rangle - \max_{v \in \partial h(x_{0})} \langle v, -d(x_{0})\rangle \\
& \stackrel{\eqref{vazio2}}{=} \langle \nabla g(x_{\alpha}), -d(x_{0})\rangle + \langle v_{1}(x_{0}), d(x_{0})\rangle \\
 & \le\big(\alpha L(x_{0}) - 1)\big)\sigma_{1}^{2}(x_{0}).
\end{aligned}
\end{equation*}
Hence $\Bar{\alpha}\geq \hat{\alpha} := 1/L(x_0)$, and therefore for all $\alpha \in [0,\hat{\alpha}]$ we get
\begin{equation*}
\xi(\alpha) = \xi(0) + \int_{0}^{\alpha} \xi'(t, 1) dt \leq \xi(0) + \int_{0}^{\alpha}(tL(x_0) - 1) \sigma_1^2(x_0) dt = \xi(0) - \frac{\alpha\big(2 - \alpha L(x_0)\big)}{2} \sigma_1^2(x_0).
\end{equation*}
By $x_{\hat{\alpha}} \in \mathcal{L}(x_0)$, the latter means that 
$$
\displaystyle \phi(x_0) - \phi(\bar{x}) \geq \phi(x_0) - \phi(x_{\hat{\alpha}}) \geq \frac{\sigma_1^2(x_0)},{2L(x_0)},
$$
which verifies \eqref{vazio4} and thus completes the proof of the lemma.
\end{proof}\vspace*{-0.05in}

Now we are ready to establish our major observation about the lower exponent PLK  property for problems of difference programming $({\cal P})$, with ${\cal C}^{1,1}$ functions $g$ and convex functions $h$, at local minimizers of $\phi$.\vspace*{0.05in}

\begin{theorem}[\bf inconsistency of the lower exponent PLK property with Lipschitz continuity of gradients]\label{lipDC}  
Let $\bar x\in{\rm int}({\rm dom}\,h)$ be a local minimizer of problem $({\cal P})$, where $g$ is of class ${\cal C}^{1,1}$ around $\bar x$, and where $h$ is convex. Then the exponent PLK property of $\phi$ at $\bar x$ fails whenever $q\in(0,1/2)$.
\end{theorem}\vspace*{-0.05in}
\begin{proof} Observe that the exponent PLK property of $\phi$ at $\bar x$ can be represented, provided that $\phi$ is locally continuous around $\bar x$ as follows from the assumptions of the theorem, in the following form: there exist a constant $M > 0$ and a neighborhood $U$ of $\bar x$ such that for all $x\in U$ taken inside of ${\cal L}(x_0)$, we have the inequality
\begin{equation}\label{vazio5}
\sigma(x) \geq M[\phi(x) - \phi(\bar{x})]^q,
\end{equation}
where $\sigma(x)$ is some estimate for the norm of the corresponding subgradient from \eqref{desigualdadeklND} of $\phi$ at $x$. The standard class $\sigma(x)$ is $\sigma(x) = \sigma_0(x)$, but we consider even a broader class where $\sigma(x) = \sigma_1(x)$ in \eqref{vazio5}. Then Lemma~\ref{lemma2} tells us that for any $x_0 \in U$, the estimates
$$ 
M \Big[\frac{1}{2 L(x_0)}\sigma_1^2(x_0)\Big]^q \leq M\big[\phi(x_0) - \phi(\bar{x})\big]^q \stackrel{\eqref{vazio5}}{\leq} \sigma_1(x_0)
$$
are satisfied by definition \eqref{exact}. This implies that
\begin{equation}\label{est}
\displaystyle M\sigma_1^{2q-1}(x_0) \leq [2L(x_0)]^{q}.
\end{equation}
It follows from \eqref{est} for $q \in (0,\frac{1}{2})$ that $L(x_0) \to \infty$ as $x_0 \to \bar{x}$,  which is inconsistent with the Lipschitz continuity  of $\nabla g$ around $\bar x$ and thus completes the proof. 
\end{proof}\vspace*{-0.05in}

As we see, the inconsistency result of Theorem~\ref{lipDC} affects also problems of {\em unconstrained minimization of ${\cal C}^{1,1}$ functions} at their local minimizers. After a preliminary version of this paper was uploaded to arXiv, we received a message from Nicolas Boumal informing us that about the results and discussions on the failure of the lower exponent PLK conditions for problems of ${\cal C}^{1,1}$-optimization given in \cite[Remark~2.21]{Boumal2023}.\vspace*{0.03in}

The following example shows that the inconsistency  observation of Theorem~\ref{lipDC}  may {\em fail} if $g$ is {\em not of class ${\cal C}^{1,1}$}. The function below is taken from \cite[Example~4.6]{li}, where it is used to demonstrate that the tight {\em quadratic} convergence rate can be achieved for the proximal methods applied to this function.\vspace*{0.05in}

\begin{example}[\bf ${\cal C}^{1,1}$ property is essential for inconsistency with lower exponent PLK] {\rm Consider the univariate function $\phi(x):=|x|^{3/2}$, which has the global minimizer $\bar x=0$. It is straightforward to compute that
\begin{equation*}
\phi'(x)=\frac{3}{2}{\rm sign}(x)|x|^{\frac{1}{2}}=\left\{\begin{array}{ll}\;\;\frac{3}{2}x^{\frac{1}{2}}&\mbox{if }\;x>0,\\
\;\;0&\mbox{if }\;x=0,\\
-\frac{3}{2}(-x)^{\frac{1}{2}}&\mbox{if }\;x<0.
\end{array}\right.
\end{equation*}
Therefore, the derivative of $\phi$ is not locally Lipschitzian around $\bar x$. On the other hand, it is not difficult to check that the lower exponent PLK property holds for this function with $\varphi(t)=t^{1-q}$ and $q=1/3$.}   
\end{example}
\vspace*{-0.1in}

\section{Conclusions}\label{sec:final}\vspace*{-0.05in}

This paper conducts a comprehensive convergence analysis of generic descent algorithms under several versions of the PLK conditions. Convergence rates are established by using the exponent PLK conditions, where some really surprising results are obtained for the case of lower exponents. It is shown, in particular, that the lower exponent PLK conditions under consideration are inconsistent  with the Lipschitz continuity of the gradient $\nabla g$ for problems of minimizing the difference functions $g-h$, where $h$ is extended-real-valued and convex. 

There are several challenging open questions for the future research. Among such questions, the following should be mentioned. First, we intend to clarify for which particular algorithms the additional assumption ensuring the striking result of Theorem~\ref{PLKfail} holds. Another issue is to investigate difference programs discussed in Section~\ref{boosted} and Section~\ref{lipDC} in the case where $h$ is nonconvex with some specific structures, e.g., being {\em prox-regular} \cite{mordukhovich2024,rw}. Among other algorithms, we plan to study from this viewpoint are the semi-Newton regularized method of \cite{aragon2023coderivative} in difference programming as well as the cubic Newton method by Nesterov and Polyak \cite{nest} and its far-going extension in \cite{li}. Let us finally mention our intention to consider DC programs where the function $g$ in the DC decomposition $\phi=g-h$ is {\em nondifferentiable}. Such problems exhibit significant differences from the differentiable case; see the recent works \cite{will,fer} for new phenomena in {\em nonexact} versions of Boosted DCA. Our goal is to establish appropriate counterparts of the results obtained in this paper for difference problems without any differentiability assumptions.\vspace*{-0.05in}

\subsection*{Acknowledgements}
Research of Boris S. Mordukhovich was supported by the US National Science Foundation under grant DMS-2204519 and by the Australian Research Council under Discovery Project DP250101112. Research of Glaydston de C. Bento was supported by CNPq grants 314106/2020-0. Research of Yurii Nesterov was supported by the National Research, Development and Innovation Office (NKFIH) under grant number 2024-1.2.3-HU-RIZONT-2024-00030.

\end{document}